\documentclass[12pt, reqno]{amsart}
\usepackage{hhline}
\usepackage{tikz}
\usepackage{verbatim}
\usepackage{amsmath}
\usepackage{amsxtra}
\usepackage{amscd}
\usepackage{amsthm}
\usepackage{amsfonts}
\usepackage{amssymb}
\usepackage{eucal}
\usetikzlibrary{arrows}
\usepackage{array}

\usepackage[all, cmtip]{xy}
\usepackage{stmaryrd}
\usepackage{color}
\usepackage{amsthm,amsfonts,amssymb,amsmath,amsxtra}
\usepackage[all]{xy}
\SelectTips{cm}{}
\usepackage{xr-hyper}
\usepackage[colorlinks=
   citecolor=Black,
   linkcolor=Red,
   urlcolor=Blue]{hyperref}
\usepackage{verbatim}

\usepackage[margin=1.25in]{geometry}

\usepackage{mathrsfs}

\RequirePackage{xspace}
\RequirePackage{etoolbox}
\RequirePackage{varwidth}
\RequirePackage{enumitem}
\RequirePackage{tensor}
\RequirePackage{mathtools}
\RequirePackage{longtable}
\RequirePackage{multirow}

\def\<{\langle}
\def\>{\rangle}

\newcommand{\sS}{\ensuremath{\mathscr{S}}\xspace}

\newcommand{\BA}{\ensuremath{\mathbb {A}}\xspace}

\newcommand{\BC}{\ensuremath{\mathbb {C}}\xspace}

\newcommand{\BF}{\ensuremath{\mathbb {F}}\xspace}
\newcommand{{\BG}}{\ensuremath{\mathbb {G}}\xspace}

\newcommand{{\BK}}{\ensuremath{\mathbb {K}}\xspace}
\newcommand{\BL}{\ensuremath{\mathbb {L}}\xspace}

\newcommand{\BP}{\ensuremath{\mathbb {P}}\xspace}
\newcommand{\BQ}{\ensuremath{\mathbb {Q}}\xspace}
\newcommand{\BR}{\ensuremath{\mathbb {R}}\xspace}
\newcommand{\BS}{\ensuremath{\mathbb {S}}\xspace}

\newcommand{\BZ}{\ensuremath{\mathbb {Z}}\xspace}

\newcommand{\CF}{\ensuremath{\mathcal {F}}\xspace}
\newcommand{\CG}{\ensuremath{\mathcal {G}}\xspace}
\newcommand{\CH}{\ensuremath{\mathcal {H}}\xspace}

\newcommand{\CM}{\ensuremath{\mathcal {M}}\xspace}

\newcommand{\CO}{\ensuremath{\mathcal {O}}\xspace}

\DeclareMathOperator{\End}{End}

\DeclareMathOperator{\Gal}{Gal}
\newcommand{\GL}{\mathrm{GL}}

\newcommand{\wh}{\widehat}

\newcommand{\ov}{\overline}

\newtheorem{theorem}{Theorem}

\theoremstyle{definition}

\newenvironment{altenumerate}
   {\begin{list}
      {\textup{(\theenumi)} }
      {\usecounter{enumi}
       \setlength{\labelwidth}{0pt}
       \setlength{\labelsep}{0pt}
       \setlength{\leftmargin}{0pt}
       \setlength{\itemsep}{\the\smallskipamount}
       \renewcommand{\theenumi}{\roman{enumi}}
      }}
   {\end{list}}
\newenvironment{altitemize}
   {\begin{list}
      {$\bullet$}
      {\setlength{\labelwidth}{0pt}
	   \setlength{\itemindent}{5pt}
       \setlength{\labelsep}{5pt}
       \setlength{\leftmargin}{0pt}
       \setlength{\itemsep}{\the\smallskipamount}
      }}
   {\end{list}}

\numberwithin{equation}{section}
\numberwithin{theorem}{section}

\setitemize[0]{leftmargin=*,itemsep=\the\smallskipamount}
\setenumerate[0]{leftmargin=*,itemsep=\the\smallskipamount}

\renewcommand{\to}{%
   \ifbool{@display}{\longrightarrow}{\rightarrow}%
   }
\let\shortmapsto\mapsto
\renewcommand{\mapsto}{%
   \ifbool{@display}{\longmapsto}{\shortmapsto}%
   }
\newlength{\olen}
\newlength{\ulen}
\newlength{\xlen}
\newcommand{\xra}[2][]{%
   \ifbool{@display}%
      {\settowidth{\olen}{$\overset{#2}{\longrightarrow}$}%
       \settowidth{\ulen}{$\underset{#1}{\longrightarrow}$}%
       \settowidth{\xlen}{$\xrightarrow[#1]{#2}$}%
       \ifdimgreater{\olen}{\xlen}%
          {\underset{#1}{\overset{#2}{\longrightarrow}}}%
          {\ifdimgreater{\ulen}{\xlen}%
             {\underset{#1}{\overset{#2}{\longrightarrow}}}
             {\xrightarrow[#1]{#2}}}}%
      {\xrightarrow[#1]{#2}}
   }
\makeatother
\newcommand{\xyra}[2][]{%
   \settowidth{\xlen}{$\xrightarrow[#1]{#2}$}%
   \ifbool{@display}%
      {\settowidth{\olen}{$\overset{#2}{\longrightarrow}$}%
       \settowidth{\ulen}{$\underset{#1}{\longrightarrow}$}%
       \ifdimgreater{\olen}{\xlen}%
          {\mathrel{\xymatrix@M=.12ex@C=3.2ex{\ar[r]^-{#2}_-{#1} &}}}%
          {\ifdimgreater{\ulen}{\xlen}%
             {\mathrel{\xymatrix@M=.12ex@C=3.2ex{\ar[r]^-{#2}_-{#1} &}}}
             {\mathrel{\xymatrix@M=.12ex@C=\the\xlen{\ar[r]^-{#2}_-{#1} &}}}}}%
      {\mathrel{\xymatrix@M=.12ex@C=\the\xlen{\ar[r]^-{#2}_-{#1} &}}}%
   }
\makeatletter
\newcommand{\xla}[2][]{%
   \ifbool{@display}%
      {\settowidth{\olen}{$\overset{#2}{\longleftarrow}$}%
       \settowidth{\ulen}{$\underset{#1}{\longleftarrow}$}%
       \settowidth{\xlen}{$\xleftarrow[#1]{#2}$}%
       \ifdimgreater{\olen}{\xlen}%
          {\underset{#1}{\overset{#2}{\longleftarrow}}}%
          {\ifdimgreater{\ulen}{\xlen}%
             {\underset{#1}{\overset{#2}{\longleftarrow}}}
             {\xleftarrow[#1]{#2}}}}%
      {\xleftarrow[#1]{#2}}
   }
\newcommand{\isoarrow}{%
   \ifbool{@display}{\overset{\sim}{\longrightarrow}}{\xrightarrow\sim}%
   }

\begin{document}

\title{The work of Peter Scholze}

\author{M. Rapoport}
\address{Mathematisches Institut der Universit\"at Bonn, Endenicher Allee 60, 53115 Bonn, Germany, and Department of Mathematics, University of Maryland, College Park, MD 20742, USA}
\email{rapoport@math.uni-bonn.de}

\date{\today}

\begin{abstract}
He has developed powerful methods in algebraic geometry over $p$-adic fields, and has proved striking theorems in this area. 
\end{abstract}

\date{\today}
\maketitle

%\tableofcontents
My purpose here is to convey some idea of the accomplishments of Peter Scholze for which he was awarded the Fields medal. 

Scholze has made ground-breaking contributions to fundamental problems in arithmetic geometry.  Although his main results so far concern the  geometry of algebraic varieties over $p$-adic fields, with important applications to the theory of automorphic forms, he has a much wider vision of mathematics. In particular, he has also contributed substantially to algebraic topology and has recently begun developing some fascinating ideas on arithmetic geometry beyond the $p$-adic setting. 

Moreover, although Scholze has made major additions to the elaborate theoretic foundations of arithmetic geometry, at the same time his ideas have dramatically simplified and clarified our field. This is a characteristic feature of his universal approach to and vision of mathematics. 

This report is structured as follows. In sections \ref{s:Frob} and \ref{s:first}, we present Scholze's perfectoid spaces and some of their first applications. In section \ref{s:proetop}, his pro-\'etale topology is introduced. This technique will be used in the proofs in sections \ref{s:Hodgerig} and \ref{s:inthodge} on $p$-adic Hodge theory and then applied, in conjunction with the Hodge-Tate period map of section \ref{s:periods}, to global problems in section \ref{s:exglob}.  Section \ref{s:vtop} is devoted to Scholze's theory of $v$-sheaves that extends the theory of diamonds from section \ref{s:proetop}.  In sections \ref{s:locshim} and  \ref{s:Lpara}, applications of these methods to  local Shimura varieties and their cohomology, and the construction of Langlands parameters are presented. 
Thus the report is organized in three themes: \emph{$p$-adic geometry} (sections \ref{s:Frob}, \ref{s:first}, \ref{s:proetop}, \ref{s:vtop}), \emph{$p$-adic Hodge theory} (sections \ref{s:Hodgerig}, \ref{s:inthodge}, \ref{s:periods}) and  \emph{(local and global) Shimura varieties and Langlands correspondences} (sections \ref{s:exglob}, \ref{s:locshim}, \ref{s:Lpara}).
Section \ref{s:further} mentions some further results of Scholze. The report ends with a short summary assessment of Scholze as a mathematician. 

\smallskip

{\tiny I thank L.~Fargues, E.~Hellmann and M.~Morrow for their remarks on this text.}

\section{Transferring the Frobenius map into mixed characteristic}\label{s:Frob}
Let $p$ be a prime number. In various aspects, algebraic  varieties in characteristic $p$, such as over $\mathbb F_p$, are easier to handle than in characteristic zero, such as over $\mathbb Q_p$. This may seem paradoxical to an analyst who works over fields of characteristic zero such as $\mathbb R$ or $\mathbb C$.  However, in characteristic $p$ the {\em Frobenius map}, mapping an element to its $p$-th power, is  compatible not only with multiplication but also with addition and therefore provides an extra symmetry which often simplifies algebraic problems. 

One of the fundamental methods developed by Scholze is his theory of \emph{perfectoid spaces}: this  presents a general framework to reduce problems about algebraic varieties in characteristic zero to algebraic varieties in characteristic $p$.

Let us give an idea of this theory. Let us start by comparing the field ${\mathbb Q}_p$ of $p$-adic numbers with the field $\mathbb F_p((t))$ of Laurent series with coefficients in the finite field $\mathbb F_p$. In the first case, elements may be written as $p$-adic expansions $\sum_ia_ip^i$, where $0\le a_i<p$, and in the second case elements may be written as $t$-adic expansions $\sum_ia_it^i$, where $a_i\in\mathbb F_p$. Thus they look superficially similar, but of course are quite different.  However, by a theorem of Fontaine-Wintenberger, after adjoining successively higher and higher $p$-power roots of  $p$, resp. $t$, these fields become more and more similar. In fact, after adjoining all the $p$-power roots, although they are not isomorphic, they have the same absolute Galois group. Scholze had the deep insight that this theorem is the manifestation of a  much more general phenomenon. 

The Fontaine-Wintenberger theorem may be reinterpreted as an equivalence between the category of finite extensions of $\mathbb Q_p(p^{1/p^\infty})$ and the corresponding category for  $\mathbb F_p((t^{1/p^\infty}))$. Scholze's perception of this theorem is that this is  merely the simplest, zero-dimensional case of a ``tilting equivalence''.  He first introduces the notion of a {\it perfectoid field}: this is a complete topological field $k$ whose topology is induced by a non-archimedean absolute value $|\,\,|\colon k\to \mathbb{R}_{\ge0}$ with dense image, such that $p$ is topologically nilpotent in $k$, and such that every element of $O_k/pO_k$ admits a $p^{\text{\scriptsize th}}$-root. Here $O_k\subseteq k$ denotes the subring of elements absolute value $\le 1$. For example, the completions of the fields $\mathbb Q_p(p^{1/p^\infty})$ and $\mathbb F_p((t^{1/p^\infty}))$ are perfectoid. Taking this as a starting point, Scholze defines a whole class of {\it perfectoid algebras} over perfectoid fields: these are certain algebras equipped with a topology, again satisfying a certain $p^{\text{\scriptsize th}}$-root condition.  And he constructs  a \emph{tilting functor} which associates to each perfectoid algebra of characteristic zero a perfectoid algebra of characteristic $p$. He shows that this is an equivalence of categories. Furthermore, he  then geometrizes this construction by introducing {\it perfectoid spaces} over a perfectoid field $k$, obtained by gluing the adic spectra of perfectoid rings (just as schemes are obtained by gluing the spectra of rings). Here the \emph{adic spectrum}, introduced by Huber in the 1990s, is  a refinement of the usual notion of spectrum in algebraic geometry which takes into account a topology on the ring. The adic spectrum of a perfectoid ring is known as an {\em affinoid perfectoid}. 

The fundamental theorem about perfectoid spaces is as follows:
\begin{theorem}\label{mainperfd}
\begin{altenumerate}
\item
Let $k$ be a perfectoid field, and denote by $k^\flat$ its tilt which is a perfectoid field of characteristic $p$. The tilting functor $X\mapsto X^\flat$ induces an equivalence of categories  between the category of perfectoid spaces over $k$ and the category of perfectoid spaces over $k^\flat$. Furthermore, the tilting functor induces an equivalence of \'etale sites, $X_{\text{{\rm \'et}}}\simeq X^\flat_{\text{{\rm \'et}}}$.  
\item
For any perfectoid space $X$, the structure pre-sheaf $\CO_X$ is a sheaf and, if $X$ is affinoid perfectoid, then 
$$
 H^i_{\text{{\rm \'et}}}(X, \CO_X)=0 \text{ for $i>0$ }.
$$
%. Furthermore, if $X$ is affinoid perfectoid,
%$$
%H^i_{\rm et}(X, \CO_X)=0, \quad \text{ for $i>0$} .
%$$
\end{altenumerate}
\end{theorem} 

The first part of the theorem is Scholze's tilting equivalence for perfectoid spaces, which simultaneously extends the Fontaine-Wintenberger theorem and the \emph{almost purity theorem} of Faltings, which was one of Faltings' key techniques in his work on $p$-adic Hodge theory. The second part of the theorem is the analogue for perfectoid spaces of Tate's acyclicity theorem for rigid-analytic spaces, or Cartan's theorem B for Stein complex spaces, or Serre's vanishing of higher cohomology on affine schemes. It is surprising in this highly non-noetherian situation.

There is also the notion of a perfectoid space without mention of a perfectoid ground field (Fontaine, Kedlaya). However, when the perfectoid ground field is not fixed, the tilting operation is not `injective': in fact, the 'moduli' of all untilts over $\BQ_p$ of a fixed complete algebraically closed field of characteristic $p$ is  the \emph{Fargues-Fontaine curve} from $p$-adic Hodge theory, a `compact $p$-adic Riemann surface' (in particular, a regular noetherian scheme of Krull dimension one) over $\BQ_p$  whose geometric properties are closely tied to $p$-adic arithmetic.

\section{First applications of perfectoid spaces}\label{s:first}
Scholze's first application of  his theory of perfectoid spaces was  a proof of Deligne's \emph{weight monodromy conjecture}  for a new class of algebraic varieties. Let $F$ be a finite extension of $\BQ_p$, and let $X$ be a proper smooth variety over $F$. Deligne's conjecture is that for any degree $i$, the monodromy filtration on the \'etale cohomology group $H^i_{\text{{\rm \'et}}}( X_{\ov F}, \BQ_\ell)$ is pure of weight $i$ (essentially, that the associated graded pieces of the monodromy filtration afford an action by the Frobenius automorphism which is pure of a certain prescribed weight).   This is undoubtedly the single most important open conjecture on the  \'etale  cohomology of algebraic varieties. Scholze proves:
\begin{theorem}
Let $X$ be a  proper smooth algebraic variety over $F$ such that $X$ is a set-theoretic complete intersection in a projective smooth toric variety. Then the weight monodromy conjecture is true for $X$. 
\end{theorem}
The proof of this theorem uses the tilting functor to reduce  subtly to  the analogous conjecture in which $F$ is replaced by a finite extension of $\mathbb F_p((t))$, which was  proved earlier by Deligne. It is conceivable that any (projective smooth) algebraic variety over $F$ satisfies the hypothesis  of Scholze's theorem, but as long as this is not known, the monodromy conjecture has to be considered as still open in the general case. 

The theory of perfectoid spaces has led to other applications. We mention a few of these in commutative algebra due to others.

\begin{altitemize}
\item The  proof of  Hochster's direct summand conjecture (Andr\'e and Bhatt).
\item The proof of  Hochster's conjecture on the existence and weak functoriality of big Cohen-Macaulay algebras (Andr\'e, Heitmann-Ma).
\item The $p$-adic analogue of Kunz's characterization of regular rings through their Frobenius endomorphism (Bhatt-Iyengar-Ma).
\item The proof of the Auslander-Goldman-Grothendieck purity conjecture on the Brauer group of schemes (\v Cesnavi\v{c}ius). 
\end{altitemize}

\section{The pro-\'etale topology and diamonds}\label{s:proetop}
One of Grothendieck's main inventions was the introduction of  the \'etale topology of schemes which lead him to a dramatic reworking of the concept of a topology. Scholze extends in several ways Grothendieck's concepts, with strong consequences. In this section we address Scholze's  pro-\'etale topology; we will do this  in the framework of perfectoid spaces,  though there are also analogues for schemes.

A morphism $f\colon {\rm Spa} (B, B^+)\to {\rm Spa} (A, A^+)$ of affinoid perfectoids is  \emph{pro-\'etale} if $(B, B^+)$ is a completed filtered colimit of perfectoid pairs $(A_i, A^+_i)$ which are \'etale over $(A, A^+)$; this definition is extended to morphisms $f\colon X\to Y$ of perfectoid spaces, so as to be local on the source and the target. In contrast to  \'etale morphisms of schemes, pro-\'etale morphisms can have infinite degree. Another subtlety  is that it may happen that the inclusion of a point in an affinoid perfectoid is  a pro-\'etale morphism. Using pro-\'etale morphisms, Scholze defines the \emph{pro-\'etale topology}. He proves the following analogue of Theorem \ref{mainperfd}. 
\begin{theorem}\label{proetdesc}
Any perfectoid space is a sheaf for the pro-\'etale topology. Furthermore, for  any perfectoid space $X$, the presheaf $\CO_X$ on the pro-\'etale site is a  sheaf  and, if $X$ is affinoid perfectoid, then 
$$
 H^i_{\text{{\rm pro-\'et}}}(X, \CO_X)=0 \text{ for $i>0$ }.
$$
\end{theorem}
%The pro-\'etale topology is a very fine topology. For instance, Scholze proves (under mild hypotheses) that any `pro-finite' morphism is pro-\'etale,  locally for the pro-\'etale topology. 
To work with the pro-\'etale topology, 
Scholze introduces the notion of a \emph{totally disconnected} perfectoid space: this is a (quasi-compact and quasi-separated)
 perfectoid space $X$ that is as close as possible to a profinite topological space, in the sense that each connected component  has a unique closed point. He proves that  any perfectoid space  may be covered, in the sense of the pro-\'etale topology, by totally disconnected ones. This  is somewhat reminiscent of the fact that any compact Hausdorff space is the continuous image of a pro-finite set. Moreover, when $X$ is totally disconnected, he  proves (roughly)  that a morphism
to $X$ is pro-\'etale if and only if its geometric fibers are profinite sets. This result gives a fiberwise criterion to decide whether  a morphism is  pro-\'etale, locally for the pro-\'etale topology on the base, and  makes the pro-\'etale topology  manageable.

The pro-\'etale topology then leads to the notion of a \emph{diamond}: a diamond is a sheaf for the pro-\'etale topology on the category of perfectoid spaces in characteristic $p$ which can be written as a quotient of a perfectoid space by a pro-\'etale equivalence relation. This definition is analogous to Artin's definition of \emph{algebraic spaces}, and expresses the intuitive idea that a diamond is obtained by glueing perfectoid spaces along pro-\'etale overlaps.  Theorem \ref{proetdesc} enables one to extend the tilting functor from perfectoid spaces to all rigid-analytic spaces: Scholze thus   defines the \emph{diamond functor} 
\begin{equation}\label{diamfunct}
\{ \text {\it  adic spaces over $\BQ_p$}\}\to \{\text{\it diamonds}\}, \quad X\mapsto X^\diamondsuit ,
\end{equation}
which, when restricted to the full subcategory of perfectoid spaces over $\BQ_p$, induces the tilting functor $X\mapsto X^\flat$. In fact, for any non-archimedean field $L$, the functor $X\mapsto X^\diamondsuit$ defines a fully faithful functor from the category of \emph{seminormal} rigid-analytic spaces over $L$ to the category of diamonds over ${\rm Spd}\, L={\rm Spa}(L, O_L)^\diamondsuit$.

The category of diamonds is much more flexible than the category of adic spaces, e.g., it allows one to take a  product of diamonds  ${\rm Spd}(\BQ_p)\times {\rm Spd}(\BQ_p)$. 
In this way, Scholze gives a meaning to the `arithmetician's dream object' ${\rm Spec} (\mathbb{Z})\times_{\BF_1} {\rm Spec} (\mathbb{Z})$ after localization at $(p,p)$, where $\BF_1$ is the non-existent field with one element.  Here the two copies of the prime number $p$ have to be thought of as two independent variables. 
 
 Scholze uses the category of diamonds also as a method to construct objects in the category of rigid-analytic spaces by first showing that these objects exist as diamonds and then showing that they are in the essential image of the diamond functor. It appears that the concept of diamonds is just the right one to address topological questions in $p$-adic geometry.

\section{Hodge theory of rigid-analytic spaces}\label{s:Hodgerig}
The classical subject of Hodge theory is concerned with the singular cohomology and de Rham cohomologies of compact complex manifolds, and their relation. It applies not only to projective algebraic varieties over $\mathbb C$ but to the wider class of compact K\"ahlerian  manifolds. The analogous {\em $p$-adic Hodge theory} of $p$-adic algebraic varieties was  initiated by Tate in the 1960s and subsequently completed by Fontaine-Messing, Faltings, Kato, Tsuji, Niziol and Beilinson.  Tate asked in his original paper whether the theory worked not only for $p$-adic varieties but for the  wider class of $p$-adic rigid-analytic spaces, which are the $p$-adic analogues of complex manifolds. The positive resolution of the main theorems of $p$-adic Hodge theory in this degree of generality is  given by the following theorem of Scholze. Here the singular cohomology groups of the classical theory are replaced by the \'etale cohomology groups.

%studying the $p$-

%\emph{$p$-adic Hodge theory} is the analogue in the $p$-adic case of the classical Hodge theory of K\"ahlerian {\bf[Capital K?]} complex manifolds. Recall that this classical theory is concerned with singular and de Rham cohomologies and their relation, not only of 

%K\"ahlerian complex manifolds, more precisely with their 

%\emph{$p$-adic Hodge theory} is the analogue in the $p$-adic case of the Hodge theory of , which was extended to general {\bf varieties} over $\BC$ by Deligne. {\bf Recall that this classical theory is concerned with singular and de Rham cohomologies, and their relation.} The Hodge theory of $p$-adic varieties {\bf [I slightly simplified ``schemes over $p$-adic schemes'']} was initiated by Tate {\bf in the 1960s} and . Tate asked more than 40 years ago whether one can transpose this theory to rigid-analytic spaces. The answer is given by the following theorem of Scholze   which is analogous to the main theorems of classical Hodge theory. 

\begin{theorem}\label{Hodgerig}
Let $X$ be a proper smooth rigid-analytic space over a complete  algebraically closed extension $C$ of $\BQ_p$.
\begin{altenumerate}
\item The Hodge-de Rham spectral sequence
$$
E_1^{i j}=H^j(X, \Omega_{X/C}^i)\Rightarrow H_{\rm dR}^{i+j}(X/C)
$$
degenerates at the first page. Moreover, for all $i\geq 0$,
$$
\sum_{j=0}^{i}\dim_C H^{i-j}(X, \Omega_{X/C}^j)=\dim_C H_{\rm dR}^i(X/C)=\dim_{\BQ_p} H_{\text{{\rm \'et}}}^i(X, \BQ_p) .
$$
\item There is a  \emph{Hodge-Tate spectral sequence}
$$
E_2^{i j}=H^i(X, \Omega_{X/C}^j)\Rightarrow H_{\text{{\rm \'et}}}^{i+j}(X, \BZ_p)\otimes_{\BZ_p} C
$$
that degenerates at the second page. 
\end{altenumerate} 
\end{theorem}

The first part of the theorem implicitly includes the statement that the \'etale cohomology groups $H_{\text{{\rm \'et}}}^i(X, \BQ_p)$ are finite-dimensional; this was for a long time conjectural. A key technique  in the proof of this theorem is that any rigid-analytic space may be covered, with respect to the pro-\'etale topology, by affinoid perfectoids. This allows one to then apply the vanishing theorems for the structure sheaves on such spaces, as in Theorem \ref{mainperfd}.

It is remarkable that, contrary to the complex case, the theorem holds  without any K\"ahler type hypothesis on $X$.

\section{Integral $p$-adic Hodge theory}\label{s:inthodge}

When the rigid-analytic space $X$ in Theorem \ref{Hodgerig} comes from a proper smooth formal scheme $\frak X$ over the ring of integers $O_{C}$, one can refine the de Rham cohomology of $X$ and prove comparison theorems with the \'etale cohomology of $X$, and also with the crystalline cohomology of the special fiber of $\frak X$. Let $A_{\rm inf}=W(O^\flat_C)$ be Fontaine's ring, with its Frobenius automorphism $\varphi$  and a fixed generator $\xi$ of $\ker(A_{\rm inf}\to O_{C})$. Also, let $k$ denote the residue field of $O_C$. 

In joint work with Bhatt and Morrow, Scholze proves: 
\begin{theorem}\label{Ainfcoho}
There exists a perfect complex $R\Gamma_{A_{\rm inf}}(\frak X)$ of $A_{\rm inf}$-modules together with a $\varphi$-linear endomorphism $\varphi\colon R\Gamma_{A_{\rm inf}}(\frak X)\to R\Gamma_{A_{\rm inf}}(\frak X)$ that becomes an automorphism after inverting $\xi$, resp. $\varphi(\xi)$. Each cohomology group $H^i_{A_{\rm inf}}(\frak X)$ is a finitely presented ${A_{\rm inf}}$-module that becomes free after inverting $p$. Furthermore, one has the following comparison isomorphisms.
\begin{altenumerate}
\item {\rm de Rham:} $R\Gamma_{A_{\rm inf}}(\frak X)\otimes^\BL_{A_{\rm inf}} O_C\simeq R\Gamma_{\rm dR}(\frak X/O_C)$. 
\item {\rm \'etale:} $R\Gamma_{A_{\rm inf}}(\frak X)\otimes_{A_{\rm inf}} W(C^\flat)\simeq R\Gamma_{\text{\rm \'et}}( X, \BZ_p)\otimes_{\BZ_p}W(C^\flat)$, $\varphi$-equivariantly.
\item {\rm crystalline:} $R\Gamma_{A_{\rm inf}}(\frak X)\otimes^\BL_{A_{\rm inf}} W(k)\simeq R\Gamma_{\rm crys}(\frak X_k/W(k))$, $\varphi$-equivariantly.
\end{altenumerate}
\end{theorem}
As a consequence of this theorem, one gets bounds for the torsion in the \'etale cohomology in terms of the crystalline cohomology:
$$
{\rm length}_{W(k)} H^i_{\rm crys}(\frak X_k/W(k))_{\rm tor}\geq {\rm length}_{\BZ_p} H^i_{\text{\rm \'et}}( X, \BZ_p)_{\rm tor} . 
$$
In particular if the crystalline cohomology is torsion free then the \'etale cohomology is torsion free as well.

The proof of this theorem uses in an essential way the Faltings almost purity theorem, cf.~section \ref{s:Frob}, enriched by a control of `junk torsion' via the Berthelot-Ogus functor $L\eta$. Contrary to the crystalline theory, in which the action of Frobenius comes from the fact that cystalline cohomology only depends on the special fiber of $\frak X$, the Frobenius action on $A_{\rm inf}$-theory is much more subtle; it ultimately comes from the Frobenius action on the tilt of $X$. 

The cohomology functor  $R\Gamma_{A_{\rm inf}}(\frak X)$ is a new cohomological invariant which cannot be obtained by a formal procedure from other previously known cohomology theories.

There is a further refinement of this result. Let $F$ be a finite extension of $\BQ_p$ contained in $C$, and assume that $\frak X$ comes by base change from a proper smooth scheme $\frak X_{O_F}$ over $O_F$. In joint work with Bhatt and Morrow, Scholze constructs a cohomology theory $R\Gamma_{\frak S}(\frak X_{O_F})$ which recovers the $A_{\rm inf}$-cohomology theory, i.e., 
$$
R\Gamma_{\frak S}(\frak X_{O_F})\otimes_{\frak S}A_{\rm inf}\simeq R\Gamma_{A_{\rm inf}}(\frak X) .
$$
Here $\frak S=W(k)[[z]]$ is the ring considered by Breuil and Kisin.   It is viewed as a subring of $A_{\rm inf}$ via the Frobenius on $W(k)$ and  by sending $z$ to the $p$-th power of a certain pseudo-uniformizer of $O_C$ (one deduced from a compatible choice of successive $p$-power roots of a fixed uniformizer of $O_F$). The proof in loc.~cit. is based on \emph{topological Hochschild homology}. That theory was given new foundations by Scholze in joint work with Nikolaus, see section \ref{s:further}, c); the flexibility of this novel version of THH theory is essential to the proof. 

Very recently, Scholze has constructed in joint work with Bhatt a new cohomology theory, \emph{prismatic cohomology}, which  clarifies the role of the Frobenius twist in the embedding of $\frak S$ into $A_{\rm inf}$ and reproves some of  the comparison isomorphisms in Theorem \ref{Ainfcoho}.

\section{Period maps}\label{s:periods}
By letting the rigid-analytic space in Theorem \ref{Hodgerig} vary, one obtains period maps. In classical Hodge theory, the trivialization of the local system defined by singular cohomology leads to a trivialization of  de Rham cohomology and hence, via the Hodge-de Rham spectral sequence, to period maps in the sense of Griffiths. The  Hodge-Tate spectral sequence of Theorem \ref{Hodgerig} leads to a new kind of period map. More precisely, Scholze proves:
\begin{theorem} Let $f\colon X\to Y$ be a proper smooth morphism of rigid-analytic spaces over a complete algebraically closed extension $C$ of $\BQ_p$. 
\begin{altenumerate}
\item Let $\BL$ be a lisse $\BZ_p$-sheaf on $X_{\rm {et}}$. Then for all $i\geq 0$, the higher direct image sheaf $R^if_*\BL$ is a lisse $\BZ_p$-sheaf on $Y_{\text{{\rm \'et}}}$. 
\item For any $i\geq 0$, there exists a perfectoid pro-\'etale cover  $\tilde Y\to Y$ such that the pull-back of  $R^if_*\BZ_p$ becomes constant; consequently, the filtration induced by the Hodge-Tate spectral sequence defines a period map of adic spaces over $C$,
$$
\pi^i_{\rm HT}\colon \tilde Y\to \CF^i_C ,
$$
where $\CF^i_C$ denotes a partial flag variety for a typical fiber of $R^if_*\BZ_p$. 
\end{altenumerate}
\end{theorem}
As an example, consider  the case of the universal elliptic curve $f\colon E\to M$ over the modular curve.  In this case, we obtain a  map $\pi_{\rm HT}\colon \tilde M\to \BP^1$ from the pro-\'etale cover $\tilde M$  trivializing $R^1f_*\BZ_p$ to the projective line.  The restriction of $\pi_{\rm HT}$ to $\BP^1\setminus \BP^1(\BQ_p)$ is a pro-finite \'etale cover, whereas the restriction to $\BP^1(\BQ_p)$ has one-dimensional fibers ($p$-adic lifts of Igusa curves). For an identification of $\tilde M$, comp.~the remark after Theorem \ref{HTshim} below. 
\\
\begin{comment}
In a more intrinsic way the preceding period map is given by a morphism from the diamond of $Y$ to a classifying pro-\'etale stack
$$
Y^\diamond \longrightarrow \ [ \mathrm{GL}_n (\mathbb{Z}_p) \backslash \CF^{i,\diamond}_C  ],
$$
where $\mathrm{GL}_n(\mathbb{Z}_p)$ is a profinite group.
The Hodge-de Rham period morphism is given himself by a morphism
$$
Y\longrightarrow \ [ \mathrm{GL}_n \backslash \CF^i_C  ]
$$
where, here, $\GL_n$ is an algebraic group. 
Contrary to Griffiths complex periods, where after trivializing the $\mathbb{Q}$-Betti cohomology one obtains a unique period map, the two fiber functors on motives over $\mathbb{Q}_p$ given by the $p$-adic \'etale cohomology and the de Rham cohomology give to two different period morphisms.

\end{comment} 
\section{Existence of global Galois representations}\label{s:exglob}
Scholze has used perfectoid methods to prove a long-standing conjecture on the construction of representations of the absolute Galois group of number  fields via the cohomology of locally symmetric spaces (conjecture of Grunewald-Ash). Let $F$ be a totally real field or a CM field. For a sufficiently small open compact subgroup $K\subset \GL_n(\BA_{F,f})$, consider the locally symmetric space
$$
X_K=\GL_n(F)\backslash [D\times \GL_n(\BA_{F,f})/K] ,
$$
where $D=\GL_n(F\otimes_{\BQ}\BR)/\BR_+K_\infty$ is the symmetric space for $\GL_n(F\otimes_\BQ\BR)$. Consider the singular cohomology groups with coefficients in $\BF_p$, for some prime number $p$.
\begin{theorem}\label{Exgal}
For any system of \emph{Hecke eigenvalues} $\psi$ appearing in $H^i(X_K, \BF_p)$, there exists a continuous semi-simple representation $\Gal(\ov F/F)\to \GL_n(\ov \BF_p)$ characterized by the property  that for all but finitely many `ramified' places $v$ of $F$, the characteristic polynomial of the Frobenius ${\rm Frob}_v$ is described in terms of the Hecke eigenvalues $\psi$  at $v$.  
\end{theorem}
In fact, a version of the theorem also holds with coefficients in $\BZ/p^m$ instead of $\BF_p$ and, passing to the limit over $m$,  yields as a consequence the existence of Galois representations in $\GL_n( \ov\BQ_p)$ attached to \emph{regular algebraic cuspidal representations} of $\GL_n(\BA_F)$ related to \emph{rational cohomology classes} proved earlier by Harris-Lan-Taylor-Thorne.   However,  rational cohomology classes are quite rare, whereas torsion classes  in the cohomology as in Theorem \ref{Exgal} abound. 

Like for that earlier result, the proof of Theorem \ref{Exgal} proceeds by realizing the cohomology of $X_K$ as the boundary contribution of a (connected) Shimura variety  of Hodge type. But by  embedding the problem into the perfectoid world, Scholze goes much farther. Let $S_K$  ($K\subset G(\BA_f)$) be a Shimura variety of Hodge type, associated to the reductive group $G$ over $\BQ$ equipped with Shimura data. Let $\{\mu\}$ be the associated  conjugacy class of cocharacters of $G_{\ov \BQ}$, and $E$  its field of definition (a finite extension of $\BQ$ contained in $\ov\BQ$).  Scholze's main tool is the following  fact. 
\begin{theorem}\label{HTshim}
Fix a prime number $p$ and a place $\frak p$ of $E$ above it. 
For any open compact subgroup $K^p\subset G(\BA_f^p)$, there exists a unique perfectoid space $S_{K^p}$ which is the \emph{completed} limit $\varprojlim_{K_p\subset G(\BQ_p)}S_{K^pK_p}\otimes_E E_{\frak p}$. Furthermore, there is a $G(\BQ_p)$-equivariant Hodge-Tate period map  $($in the sense of section \ref{s:periods}$)$, 
$$
\pi_{\rm HT}\colon S_{K^p}\to \CF_{G, \{\mu\}}\otimes_E E_{\frak p} .
$$
\end{theorem}  
In the case of the modular curve we have $G=\GL_2$ and $E=\BQ$. Then we obtain the Hodge-Tate period map mentioned at the end of section \ref{s:periods}. 

\begin{comment}
In fact one of the main technical difficulties in the proof of theorem \ref{Exgal} is the proof by Scholze that \emph{$\pi_{HT}$ extends to the minimal compactification} of the Shimura variety. Moreover the proof of \ref{Exgal} is purely $(p,p)$ (the space is $p$-adic and the coefficients are $p$-torsion), and uses the full force of his $p$-adic comparison theorems (like theorem \ref{s:Hodgerig}).
\end{comment}

As an application of these methods, Scholze also proves the following vanishing theorem,  conjectured by Calegari and Emerton.  
Recall the definition of the compactly supported \emph{completed cohomology} groups for a fixed tame level $K^p\subset G(\BA_f^p)$,
$$
\tilde H^i_{c}(S_{K^p}, \BZ_p):=\varprojlim_m\varinjlim_{K_p} H^i_c(S_{K^p K_p}, \BZ/p^m) .
$$
\begin{theorem}
For $i>\dim S_K$, the completed cohomology group with compact supports $\tilde H^i_{c}(S_{K^p}, \BZ_p)$ vanishes. 
\end{theorem}
Even without passing to the limit, one has a vanishing theorem, proved by Scholze in joint work with Caraiani: 
\begin{theorem}
Let $S_K$ be a \emph{simple} Shimura variety associated to a fake unitary group $($then  $S_K$ is  compact$)$. Let $\ell\neq p$. The localization $H^i(S_K, \mathbb{F}_\ell)_{\frak m}$ at a $p$-\emph{generic maximal ideal} of the Hecke algebra  vanishes for $i\neq \dim S_K$. 
\end{theorem}

 This result is a torsion analog of a well-known archimedean result that states that  automorphic representations $\Pi$, contributing to the singular cohomology of $S_K$, with  tempered archimedean component, only show up in the middle degree.

\begin{comment}
The main tool in the proof of the preceding theorem is a new method introduced by Scholze to study the cohomology of Shimura varieties.
Here, contrary to before, Scholze works in a $(p,\ell)$ situation (the space is $p$-adic and the coefficients are $\ell$-torsion). 
 Classically, like in the work of Harris and Taylor, this can be done by studying the complex of nearby cycles $R\psi \mathbb{F}_\ell$, that is to say '$R sp_*\mathbb{F}_\ell$', $sp$ being the specialization morphism 
with respect to the choice of an integral model over the $p$-adic numbers.  Scholze replaces this specialization morphism by $\pi_{HT}$  and studies the complex $$R\pi_{HT *}\mathbb{F}_\ell.$$ In particular he proves that this is \emph{perverse} with respect to a Newton stratification on the Hodge-Tate flag variety, and computes its fibers. The advantage of this method is that this \emph{does not involve any choice of an ad-hoc integral model}, the study is done in characteristic $0$ ! This promising approach should lead to plenty new results on the torsion in the cohomology of Shimura varieties.

\end{comment}

\section{The $v$-topology and \'etale cohomology of diamonds}\label{s:vtop}

Scholze introduces another topology on the category of perfectoid spaces, besides the pro-\'etale topology of section \ref{s:proetop}.  The  \emph{$v$-topology} is the topology obtained by declaring that any surjective map between  affinoid perfectoids is an open cover. Even though it may appear at first sight that the $v$-topology admits far too many open covers to be  useful, Scholze uses this topology  to dramatic effect: in particular, it allows him to extend the diamond functor from rigid-analytic spaces to formal schemes.  The basis of all applications  is a descent theorem for the $v$-topology:
\begin{theorem}\label{vdesc}
\begin{altenumerate}
 \item Any diamond satisfies the sheaf axioms for the  $v$-topology.  
 \item For  any perfectoid space $X$ the presheaf $\CO_X$ on the $v$-site is a sheaf,  and if $X$ is affinoid perfectoid, then 
$$
 H^i_{v}(X, \CO_X)=0 \text{ for $i>0$ }.
$$
\item For any perfectoid space $X$, the category of locally  free $\CO_X$-modules of finite rank satisfies descent for the $v$-topology. The same holds for the category of  separated \'etale morphisms.
\end{altenumerate}
\end{theorem} The statement (i) is the analogue of Gabber's theorem that any algebraic space is a fpqc-sheaf. The statement (iii) is a key tool in the work of Fargues-Scholze on the $v$-stack of vector bundles on the Fargues-Fontaine curve and its  \'etale cohomology, comp. section \ref{s:Lpara}.

Scholze also shows, under certain hypotheses,  that any  $v$-sheaf which is suitably covered  by a perfectoid space  is automatically a diamond. This is the analogue of Artin's theorem on algebraic spaces, reducing smooth, and even flat, groupoids to \'etale groupoids.

Using these concepts, Scholze has established an \emph{\'etale cohomology theory of diamonds}, taking as a model Grothendieck's \'etale cohomology theory for schemes. In particular, he  constructs the analogue of the `six-operation calculus'  and appropriate versions of the proper and smooth base change theorems. This theory is one of the key  tools in the geometric construction of smooth representations of $p$-adic groups and in the geometric construction of Langlands parameters, cf.~section \ref{s:Lpara}.  Remarkably,  for perfectoid spaces the notion of smoothness is highly non-obvious (the usual characterizations, via differentials or via infinitesimal liftings, lose their  sense in this context).

\section{Local Shimura varieties}\label{s:locshim}
A \emph{local Shimura datum} is a triple $(G, \{\mu\}, b)$ consisting of a reductive group $G$ over $\BQ_p$, a conjugacy class $\{\mu\}$ of minuscule cocharacters of $G_{\ov \BQ_p}$, and an element $b\in G(\breve \BQ_p)$ whose $\sigma$-conjugacy class lies in $B(G, \{\mu\})$, i.e.~is \emph{neutral acceptable} wrt.~$\{\mu\}$.  Here, for any finite extension $F$ of $\BQ_p$ contained in $\ov\BQ_p$, we denote by $\breve F$ be the completion of the maximal unramified extension of $F$. Let $E$ be the field of definition of $\{\mu\}$, a finite extension of $\BQ_p$ contained in the fixed algebraic closure $\ov\BQ_p$. Partly in joint work with Weinstein, Scholze proves:
\begin{theorem}
There exists a \emph{local Shimura variety} associated to $(G, \{\mu\}, b)$: a tower $\CM_{(G, \{\mu\}, b), K}$ of rigid-analytic spaces over $\breve E$, parametrized by open compact subgroups $K\subset G(\BQ_p)$, equipped with \'etale covering maps
$$
\CM_{(G, \{\mu\}, b), K}\to \breve{\CF}^{\rm adm}_{G, \{\mu\}}\subset  \breve{\CF}_{G, \{\mu\}} ,
$$
with geometric fibers $G(\BQ_p)/K$. 
\end{theorem}
Here $\breve{\CF}_{G, \{\mu\}}$ denotes the partial flag variety over $\breve E$ associated to $G$ and $\{\mu\}$, and $\breve{\CF}^{\rm adm}_{G, \{\mu\}}$ denotes the open adic subset of \emph{admissible points}. In the case that $K$ is a parahoric subgroup, Scholze constructs  a natural integral model over $O_{\breve E}$ of  $\CM_{(G, \{\mu\}, b), K}$ as a $v$-sheaf. 

The proof of this theorem proceeds by first constructing the diamond over $\breve E$ associated to $\CM_{(G, \{\mu\}, b), K}$ and then showing that it lies in the image of the fully faithful functor \eqref{diamfunct}.  The diamond is the moduli space of \emph{$p$-adic shtukas}, the $p$-adic analogue of Drinfeld's shtukas in the function field case (except that here there is  only one \emph{leg}).  Examples of  local Shimura varieties  are given  by Rapoport-Zink moduli spaces of $p$-divisible groups inside a given  quasi-isogeny class (and their integral models for parahoric level exist in this case as \emph{formal schemes} and not merely as $v$-sheaves).  This fact is highly non-trivial and 
is based on the following  description of  $p$-divisible groups due to Scholze and Weinstein which is reminiscent of Riemann's description of complex tori:
\begin{theorem}
Let $C$ be an algebraically closed complete extension of $\BQ_p$, and $O_C$ its ring of integers. There is an equivalence of categories
$$
\begin{aligned}
\{\text{$p$-divisible groups over $O_C$}\}\simeq &\{\text{pairs $(\Lambda, W)$, where $\Lambda$ is a finite free $\BZ_p$-module}\\
&\text{\quad and $W\subset \Lambda\otimes C$ is a $C$-subvector space} \} .
\end{aligned}
$$
\end{theorem}
This description of $p$-divisible groups over $O_C$ is closely related to Fargues' earlier description  in terms of \emph{integral $p$-adic Hodge theory} in the sense of section \ref{s:inthodge}.

This new point of view  of Rapoport-Zink spaces allows Scholze to establish isomorphisms between various such spaces (and their inverse limits over shrinking $K$) that have been conjectured for a long time:
\begin{altitemize}
\item $\CM_{(G, \{\mu\}, b), \infty}\simeq \CM_{(G^\vee, \{\mu^\vee\}, b^\vee), \infty}$, where $(G^\vee, \{\mu^\vee\}, b^\vee)$ denotes the \emph{dual local Shimura datum}, provided that $b$ is \emph{basic}. This solves a conjecture of Gross and Rapoport-Zink.   The case $G=\GL_n$ was proved earlier by Faltings (and Fargues), but in a more complicated indirect formulation. Furthermore, this \emph{duality isomorphism} exchanges the Hodge-Tate period map with the de Rham period map, cf.~section \ref{s:periods} (here the de Rham cohomology is trivialized). 
\item identification,  in the `fake' Drinfeld case, of a connected component of $\CM_{(G, \{\mu\}, b), K}$ with Drinfeld's formal halfspace $\wh{\Omega}^n_F\wh\otimes_F\breve F$. Here the integral $p$-adic Hodge theory in the sense of  section \ref{s:inthodge} plays a key role. 
\end{altitemize}

\section{The cohomology of local Shimura varieties and smooth representations}\label{s:Lpara}
Let $I$ be a finite set with $m$ elements. A \emph{local Shtuka datum} with $m$ legs  is a triple $(G, \{\mu_i\}_i, b)$ consisting of a reductive group $G$ over $\BQ_p$, a  collection $\{\mu_i\}_{i\in I}$ of cocharacters of $G_{\ov\BQ_p}$ and $b\in G(\breve{\BQ}_p)$. When $I=\{*\}$ and $\{\mu\}=\{\mu_*\}$ is minuscule, one recovers the definition of a local Shimura datum, cf. last section.  Generalizing the case of local Shimura varieties, Scholze constructs a tower of  diamonds (for varying $K\subset G(\BQ_p)$), 
$$f_K\colon \CM_{({G, \{\mu_i\}, b)}, K}\to\prod\nolimits_{i\in I}{\rm Spd}\, \breve E_i ,
$$ which is a moduli space of \emph{shtukas with $m$ legs bounded by $\{\mu_i\}$}.  Let $J_b$ be the $\sigma$-centralizer group of $b$ (an inner form of $G$ over $F$ when $b$ is basic). Then $J_b(\BQ_p)$ acts on each member $ \CM_{({G, \{\mu_i\}, b)}, K}$ of the tower, whereas $G(\BQ_p)$ acts on the tower as a group of Hecke correspondences. 

The tower  is equipped with a \emph{period map} to a Schubert variety inside a version of the Beilinson-Drinfeld   affine Grassmannian. When $m=1$, this Beilinson-Drinfeld Grassmannian  can be identified with the $B^+_{\rm dR}$-Grassmannian 
 of Scholze, with point set $G\big(B_{\rm dR}(C)\big)/G\big(B_{\rm dR}^+(C)\big)$ over a complete algebraically closed extension $C$ of $\BQ_p$. Here $B_{\rm dR}(C)$ and $B^+_{\rm dR}(C)$ are Fontaine's rings associated to $C$. 

Let $\Lambda$ be the ring of integers in a finite extension of $\BQ_\ell$. In their recent joint work, Fargues and Scholze associate to $\{\mu_i\}_{i\in I}$ a sheaf of $\Lambda$-modules $\sS_{\{\mu_i\}}$ on $\CM_{({G, \{\mu_i\}, b)}, K}$, to which the $J_b(\BQ_p)$-action is lifted. This construction uses the period map mentioned above. When $I=\{*\}$ and  $\{\mu\}=\{\mu_*\}$ is minuscule, then $\sS_{\{\mu_*\}}=\Lambda$. They prove the following fundamental finiteness theorem. 
\begin{theorem}\label{finiteL}
\begin{altenumerate}
\item The complex  $Rf_{K !}\sS_{\{\mu_i\}}$ comes in a natural way from an object of $D(J_b(\BQ_p)\times \prod_{i\in I} W_{E_i}, \Lambda)$, and its restriction to $D(J_b(\BQ_p), \Lambda)$ is compact $($i.e., lies in the thick triangulated subcategory generated by the $\ell$-adic completions of $\text{c-${\rm Ind}_K^{J_b(\BQ_p)} \Lambda$}$ as $K$ runs through open pro-$p$-subgroups of $J_b(\BQ_p)$$)$.
\item 
Let $\rho$ be an admissible smooth representation of $J_b(\BQ_p)$ with coefficients in $\Lambda$. Then  the object 
$$
{\rm RHom}_{J_b(\BQ_p)}(Rf_{K !}\sS_{\{\mu_i\}}, \rho)
$$
of $D(\prod_{i\in I}W_{E_i}, \Lambda)$ is a representation of $\prod_{i\in I}W_{E_i}$ on a perfect complex of $\Lambda$-modules. 
\item Passing to the limit over $K$, 
$$
\varinjlim_K\,{\rm RHom}_{J_b(\BQ_p)}(Rf_{K !}\sS_{\{\mu_i\}}, \rho)
$$
gives rise to a complex of admissible $G(\BQ_p)$-representations equipped with an action of $\prod_{i\in I}W_{E_i}$. If $\rho$ is a compact object of $D(J_b(\BQ_p), \Lambda)$, then so is this last complex of $G(\BQ_p)$-representations.
\end{altenumerate} 
\end{theorem}

One application of Theorem \ref{finiteL} is due to   Fargues and Scholze and concerns  \emph{local $L$-parameters}. This application is inspired by the work of V.~Lafforgue in the global function field case.  Let us sketch it. 

Let $I$ be a finite set, and let $V\in {\rm Rep}_\Lambda ((^L G)^I)$. Fargues and Scholze  construct a variant $\CM_{({G, V, 1)}, K}$ of $\CM_{({G, \{\mu_i\}, b)}, K}$ for $b=1$ (then $J_b=G$): a space of shtukas  bounded by $V$,
$$
f_K\colon \CM_{({G, V, 1)}, K}\to ({\rm Spd}\, \breve \BQ_p)^I , 
$$
 which is equipped with a  version of the period map.
Furthermore, Fargues and Scholze construct a sheaf of $\Lambda$-modules $\sS_V$ on $\CM_{({G, V, 1)}, K}$, to which the action of $G(\BQ_p)$ is lifted. Restriction to the diagonal 
$$
\Delta\colon {\rm Spd}\, \breve \BQ_p\to ({\rm Spd}\, \breve \BQ_p)^I 
$$
yields a moduli space of shtukas  with one leg, $f_K^\Delta\colon \CM_{({G, \Delta^*V, 1)}, K}\to{\rm Spd}\, \breve \BQ_p$, with a sheaf $\sS_{\Delta^*V}$. 

Let $i\colon G(\BQ_p)/K=\CM_{({G, \Lambda, 1)}, K}\hookrightarrow \CM_{({G, \Delta^*V, 1)}, K}$ be the subspace of shtukas with no legs.  Let $\alpha\colon\Lambda\to \Delta^* V$ and $\beta\colon \Delta^* V\to \Lambda$ be  maps of $^LG$-modules. Then $\alpha$, resp. $\beta$, induce  maps $\alpha\colon i_*\Lambda\to \sS_{\Delta^* V}$, resp.  $\beta\colon  \sS_{\Delta^* V}\to i_*\Lambda$. Let $(\gamma_i)_{i\in I}\in W_{\BQ_p}^I$, and let  $\ov x$ be a geometric point of $\Delta({\rm Spd}\, \breve \BQ_p)$.  Then we obtain the   endomorphism 
\begin{equation*}
\begin{aligned}
\text{$c$-${\rm Ind}_K^{G(\BQ_p)} \Lambda\xrightarrow{\alpha} (Rf_{K !}^\Delta\sS_{\Delta^* V})_{\ov x}=$}&\text{$(Rf_{K !}\sS_{ V})_{\ov x}$}\xrightarrow{(\gamma_i)} \\\text{$\to (Rf_{K !}\sS_{ V})_{\ov x}=$}&\text{$(Rf_{K !}^\Delta\sS_{\Delta^* V})_{\ov x}\xrightarrow{\beta}$ }\text{$c$-${\rm Ind}_K^{G(\BQ_p)} \Lambda$ .}
\end{aligned}
\end{equation*}
Here the action of $(\gamma_i)$ is given by Theorem \ref{finiteL}. Fargues and Scholze prove that this endomorphism is given by a central element of the Hecke algebra $\CH(G, K)=\End_{G(\BQ_p)}(\text{$c$-${\rm Ind}_K^{G(\BQ_p)}\Lambda)$}$. Passing to the limit over all $K$, they  define thus  an element of the Bernstein center of $G(\BQ_p)$. 

The following theorem associates  $L$-parameters to smooth representations of $G(\BQ_p)$. 
\begin{theorem}
For any irreducible smooth $\ov\BQ_\ell$-representation $\pi$ of $G(\BQ_p)$ which admits an invariant $\ov\BZ_\ell$-lattice, there is a unique  $($up to conjugation by $G^\vee(\ov\BQ_\ell)$$)$ continuous semisimple map
$$
\phi_\pi\colon W_{\BQ_p}\to ^L\! G(\ov\BQ_\ell) ,
$$
compatible with the projection of  $^L G(\ov\BQ_\ell)$ to $W_{\BQ_p}$, with the following property. Let  $I$ be a finite set, let  $V\in {\rm Rep}_\Lambda ((^L G)^I)$, let $\alpha\colon\Lambda\to \Delta^* V$ and $\beta\colon \Delta^* V\to \Lambda$, and let $(\gamma_i)_{i\in I}\in W_{\BQ_p}^I$. Then the action on  $\pi$ of the element of the Bernstein center above  is given by the scalar that arises as the composition
$$
\ov\BQ_\ell\xrightarrow{\alpha} \Delta^*V=V\xrightarrow{(\phi_\pi(\gamma_i))} V=\Delta^*V\xrightarrow{\beta}\ov\BQ_\ell .
$$
\end{theorem}
Much work remains to be done to better understand  this construction of Fargues-Scholze: which $L$-parameters arise in this way, how  the corresponding $L$-packets are related to those known by explicit representation-theoretic methods, how  these $L$-parameters can be used to construct stable distributions, etc. 

Another application of Theorem \ref{finiteL} is due to Kaletha and Weinstein, and concerns the Kottwitz conjecture on the cohomology of local Shimura varieties. Let $(G, \{\mu\}, b)$ be a local Shimura datum  such that $b$ is \emph{basic}, and let $\CM_{(G, \{\mu\}, b), K}$ be the associated local Shimura variety, cf.~section \ref{s:locshim}. As mentioned above, the Fargues-Scholze sheaf $\sS_{\{\mu\}}$ on $\CM_{(G, \{\mu\}, b), K}$ is in this case  the constant sheaf $\Lambda$. Let $\rho$ be an irreducible admissible smooth representation of $J_b(\BQ_p)$ with coefficients in $\Lambda$, and form the virtual smooth representation of $G(\BQ_p)$ from  Theorem \ref{finiteL}, (iii),
$$
H^*(G, \{\mu\}, b)[\rho]=\sum\nolimits_i (-1)^i \big(\varinjlim_K {\rm Ext}^i_{J_b(F)}(Rf_{K !}\Lambda, \rho)\big) .
$$ 
In the following theorem, $d$ denotes the dimension of $\CM_{(G, \{\mu\}, b), K}$. 
\begin{theorem}[Kaletha, Weinstein]
Let $\phi\colon W_F\to ^L\!\!G$ be a \emph{discrete Langlands parameter} for $G$, and let $\rho\in\Pi_\phi(J_b)$. Then the following identity holds in the quotient ${\rm Groth}(G(\BQ_p))^{\rm ell}$ of the Grothendieck group of admissible smooth representations of $G(\BQ_p)$ on $\ov\BQ_\ell$-vector spaces by the subgroup generated by \emph{non-elliptic} representations,
$$
H^*(G, \{\mu\}, b)[\rho]=(-1)^d\sum_{\pi\in\Pi_\phi(G)}\delta(\pi, \rho)\pi . 
$$
\end{theorem}
 Here $\Pi_\phi(G)$, resp. $\Pi_\phi(J_b)$, denotes the $L$-packet associated to the parameter $\phi$, assuming that this concept is defined (and satisfies some natural properties), as  e.g.,  for $G={\rm GSp}_4$. There is good hope that there is a definition for any tamely ramified group $G$ and $p$ sufficiently  large wrt. $G$. Whether the Fargues-Scholze definition of $\Pi_\phi(G)$, sketched above, can be used is an open problem.
 
The multiplicities ocurring here were defined earlier by Kottwitz under more restrictive hypotheses. This theorem is the confirmation of a conjecture of Kottwitz in a weaker form (weaker because the Weil group action is disregarded, and because only the image in ${\rm Groth}(G(F))^{\rm ell}$ is considered). There is also an extension of this theorem to the case when $\{\mu\}$ is no longer  assumed to be minuscule: instead of a local Shimura variety one uses the moduli space of shtukas with one leg (a diamond, but no longer a rigid-analytic space).  This uses the \emph{geometric Satake equivalence for the $B^+_{\rm dR}$-Grassmannian} of Fargues and Scholze.

\section{Further achievements}\label{s:further}

$a)$ One of Scholze's first accomplishments was a new proof of the \emph{local Langlands conjecture} for $\GL_n$ over a $p$-adic field $F$. This conjecture states that continuous representations (on finite-dimensional $\BC$-vector spaces) of the absolute Galois group $\Gal(\ov F/F)$ correspond to irreducible admissible representations of $\GL_n(F)$. It was first proved by Harris-Taylor and by Henniart. These earlier  proofs are similar and use in an essential way Henniart's \emph{numerical local Langlands conjecture}; this result in turn is based on a complicated reduction modulo $p$ method and relies ultimately on Laumon's results on the Fourier-Deligne transform and Kazhdan's construction of exotic $\ell$-adic Galois representations for function fields. Scholze's proof is purely in characteristic zero and structurally much simpler,  based instead on a geometric argument via the \emph{nearby cycles sheaves} of certain moduli spaces of $p$-divisible groups. 

\smallskip

$b)$ Let $G$ denote a reductive group over the Laurent series field $k((t))$, where $k$ is a field. Then for any parahoric group scheme $\CG$ over $k[[t]]$, one has the construction of the associated partial affine flag variety $\CF_\CG=LG/L^+\CG$, an ind-projective ind-scheme over $k$. In the case $G=\GL \big(k((t))^n\big), \CG=\GL \big(k[[t]]^n\big)$, this yields the \emph{affine Grassmannian} ${\rm Gr}^{\rm aff}=\CF_\CG$ which  parametrizes $k[[t]]$-lattices in $k((t))^n$. X.~Zhu has transposed this ``equal characteristic'' theory to the unequal characteristic: he constructs a {\emph{Witt vector Grassmannian} ${\rm Gr}^{W,{\rm aff}}$ which is an inductive limit of perfections of algebraic spaces and whose $R$-valued points, for  perfect rings $R$ of characteristic $p$, parametrize $W(R)$-lattices in $\big(W(R)[1/p]\big)^n$. Scholze, in joint work with Bhatt, shows that ${\rm Gr}^{W,{\rm aff}}$ is an ind-projective scheme, by constructing an analogue of the natural ample line bundle on ${\rm Gr}^{\rm aff}$. The main tool in this construction is the $v$-descent of vector bundles, cf.~Theorem \ref{vdesc}, (iii). Scholze interprets ${\rm Gr}^{W,{\rm aff}}$ as the special fiber of an integral model of his $B^+_{\rm dR}$-Grassmannian, cf.~section \ref{s:Lpara}.
\smallskip

$c)$ \emph{Cyclic homology} was introduced  in the early eighties to serve as an extension of de Rham cohomology to a non-commutative setting. It relies on the \emph{algebraic theory} of Hochschild homology. \emph{Topological Hochschild homology} (THH) is Hochschild homology \emph{relative to the sphere spectrum $\BS$}.   In joint work with Nikolaus, Scholze gives a  definition of \emph{Topological Cyclic homology} (TCH) in terms of a Frobenius operator on THH. This approach  avoids the ad hoc methods used earlier to define TCH,   by staying strictly within the realm of homotopy theory.  In particular, it constructs a Frobenius map in stable
homotopy theory that lives inherently in mixed characteristic, whereas the
classical Frobenius map is restricted to characteristic $p$. The relevance for algebraic geometry is furnished by the work of Bhatt-Morrow-Scholze which defines Òmotivic filtrationsÓ on THH and related theories, and relates the graded pieces with  $p$-adic cohomology theories such as crystalline cohomology and the $A_{\rm inf}$-cohomology, cf.~section \ref{s:inthodge}.   

\smallskip

\smallskip

$d)$ The existence of  a \emph{$p$-adic local Langlands correspondence} for $\GL_n(F)$ was envisioned by Breuil and was established by Colmez, Pa{\v{s}}k{\=u}nas and others in the case of $\GL_2(\BQ_p)$. Starting with a $p$-adic representation  of $\GL_n(F)$, for any $n$ and any finite extension $F$ of $\BQ_p$, Scholze produces in a purely local way a $p$-adic Galois representation. He shows that this indeed generalizes the earlier construction for $n=2$ and $F=\BQ_p$, and also relates in the latter  case this local construction to a global construction (\emph{local-global compatibility}). Scholze's proof is based on the perfectoid space associated to the limit of the Lubin-Tate tower, and the crystalline period map to $\breve{\BP}^{n-1}$, cf.~section \ref{s:locshim}. Much work remains to better understand Scholze's construction. 

\smallskip

$e)$ I refer to Scholze's write-up of his plenary lecture at this congress \cite{Sch-ICM2018} for his recent  ideas which go far beyond the $p$-adic world for a fixed $p$.

\section{summary}
Scholze has proved a whole array of theorems in $p$-adic geometry. These theorems are not disjoint but, rather, are the outflow of a theoretical edifice that Scholze has created in the last few years.

\begin{comment}
Very often, mathematicians are divided into two classes: the theory builders and the problem solvers. Of course, there is no sharp dividing line between the two: theory builders may well use their theory to solve problems and problem solvers may develop methods which are at the origin of new  theories.  Still, the distinction does exist and does make sense. Scholze is clearly a theory builder---but one with a keen sense for the striking theorem.  

\end{comment}
There is no doubt that Scholze's ideas will keep mathematicians busy for many years to come. What is remarkable about Scholze's approach to mathematics is the ultimate simplicity  of his ideas. Even though the execution of these ideas demands great  technical power (of which Scholze has an extraordinary command), it is still true that the initial key idea and the final result have the appeal of inevitability of the classics, and their elegance. We surely can expect  more great things of Scholze in the future, and it will be fascinating to see to what further heights Scholze's   work  will take him.

 \end{document}